\documentclass[11pt]{article}
\usepackage{amssymb}
\usepackage{latexsym}

\newtheorem{thm}{Theorem}[section]
\newtheorem{prop}[thm]{Proposition}
\newtheorem{lemma}[thm]{Lemma}

\newtheorem{dfn}[thm]{Definition}

\newtheorem{rmk}[thm]{Remark}
\newtheorem{ex}[thm]{Example}

\newcommand{\reals}{\mathbb R}

\newcommand{\call}{{\cal L}}

\newcommand{\mathbold}{\bf}

\newcommand{\G}{{\Lambda}}

\def\qed{\rule{2.3mm}{2.3mm}}

\newcommand{\lag}{\langle}
\newcommand{\rg}{\rangle}

\begin{document}

\title{\bf Conformal Dirac Structures }

\author{ A\"{\i}ssa WADE 
%\footnote{} 
 }

\date{}

\maketitle

\begin{abstract}
 The Courant bracket defined originally on the sections of the
 vector bundle $TM \oplus T^*M \rightarrow M $ is extended to the direct sum
  of the 1-jet vector bundle and its dual. The extended bracket allows
 one to interpret many structures encountered in differential geometry
 in terms of Dirac structures. We give here a new approach to conformal
 Jacobi structures. 

\end{abstract}

\section{Introduction}
 
  In this paper, we describe globally
  Jacobi manifolds through the theory of Dirac structures.
 Dirac structures on manifolds, introduced by T. Courant and A. Weinstein
 in~\cite{CW}, provide a framework for the study of Poisson structures,  
 pre-symplectic manifolds, as well as foliations. Roughly speaking, 
 a {\sl Dirac structure} on a smooth manifold $M$ is a sub-bundle 
   $L \rightarrow M$ of the vector bundle 
 $TM \oplus T^*M \rightarrow M$, whose total space is maximally
 isotropic under the canonical symmetric 2-form on  $TM \oplus T^*M$ and
  whose smooth sections are closed under the {\sl Courant bracket} 
  (see Section 2 below). On the other hand, a {\sl Jacobi structure} on
  $M$ is defined  by a bivector field $\pi$ and a vector field $E$ 
 satisfying $[\pi,\pi]_s=2E \wedge \pi$,  $[E,\pi]_s=0, $ 
  where $[.,.]_s$ is the Schouten-Nijenhuis bracket (~\cite{Kz}).
 A {\sl Jacobi manifold} is a smooth manifold endowed with a 
 Jacobi structure.
 It is known that Jacobi structures include contact forms,
  symplectic and Poisson  structures (see for instance~\cite{GL}). Thus a  
 natural question arises: is there a simple characterization of Jacobi 
manifolds by means of Dirac structures?\par

 The main problem which emerges when one tries 
 to describe a general Jacobi structure using the concept of a 
 Dirac structure on a manifold comes from the fact that a Jacobi structure
 on a smooth manifold $M$ involves the vector  bundle of the  1-jets of 
 real functions on  $M$ instead of the tangent bundle. 
  Hence the first step in resolving the above question is to find a suitable 
 bracket on the sections of the vector bundle 
 ${\cal E}^1(M)=(TM \times \reals) \oplus (T^*M \times \reals) \rightarrow M$ 
  extending the Courant bracket. In Section 3, we introduce 
 such a bracket on sections of ${\cal E}^1(M)$, which turns out to be an 
 operation leading to the 1-jet Lie algebroid structure associated to any 
 Jacobi manifold (see~\cite{KS}). Our approach is based on an  
 idea that can be found in~\cite{D}, where the author
  generalizes the Courant bracket to the case of complexes over Lie
 algebras, which is the algebraic counterpart of the de Rham complex of 
  differential forms. \par

   Our principal aim is to develop the theory of 
  conformal Jacobi structures (see for instance~\cite{DLM}) within the context
 of  Dirac structures. So we expect to find an equivalence
 relation for Dirac structures generalizing the concept of
 a conformal Jacobi manifold. Our motivation for this investigation originates
 from the fact that an important class of manifolds, namely contact
 manifolds without a globally defined contact 1-form, are conformal Jacobi
 manifolds in the sense of~\cite{DLM}  but not Jacobi manifolds.  \par
  Our presentation goes as follows:  we 
 introduce basic definitions and give some examples
 of Dirac structures on manifolds in Section 2. \par
  In Section 3, we consider a skew symmetric bracket  $[.,.]$ on the sections 
 of ${\cal E}^1(M)$. This bracket is just a straightforward extension 
 of the Courant bracket. 
  Next, we consider ${\cal E}^1(M)$ endowed with its canonical symmetric
 2-form and we define a Dirac structure of ${\cal E}^1(M)$ to be a maximal 
 isotropic sub-bundle $L$ of ${\cal E}^1(M)$ whose sections  
 are closed under this bracket $[.,.]$. We prove in Theorem~\ref{algebroid}
  that any Dirac structure of ${\cal E}^1(M)$ is a Lie algebroid. \par 
 In Section 4,  examples of Dirac 
 structures of ${\cal E}^1(M)$, such as homogeneous Poisson manifolds,
  Jacobi structures and Nambu manifolds are presented. \par
Section 5 is devoted to the discussion of admissible functions of a
 Dirac structure of ${\cal E}^1(M)$. We show that there is an one-to-one
  correspondence between Jacobi structures on a manifold $M$ and 
 Dirac structures of  ${\cal E}^1(M)$ for which any smooth function on  
 $M$ is admissible. \par 
 In Section 6,  we introduce an equivalence relation
  among Dirac structures of ${\cal E}^1(M)$ and  we use it
 to define a conformal Dirac structure. This section contains our 
 main results (Theorem~\ref{step2} and Theorem~\ref{Laconf}).

\section{Courant Bracket - Dirac Structures}

Let $M$ be a smooth manifold. Consider the vector bundle $TM \oplus T^*M \rightarrow M$.  The {\sl Courant bracket} on the space $\Gamma(TM \oplus T^*M)$  of smooth sections of $TM \oplus T^*M$ is the skew symmetric operation $[.,.]_c$ such that  for any $X_1+ \alpha_1$,  $X_2+ \alpha_2 \in
 \Gamma(TM \oplus T^*M)$, we have

$$[X_1+ \alpha_1, X_2+ \alpha_2]_c=[X_1,X_2]+ \call_{X_1}\alpha_2 - \call_{X_2}\alpha_1
+{1 \over 2}d(i_{X_2} \alpha_1 - i_{X_1} \alpha_2),$$

\noindent  where $i_X$ is the interior product by the vector field $X$, 
 $d$ is the exterior derivative and $\call_X=d \circ i_X + i_X \circ d$ 
  is the Lie derivation by $X$. 
 The Courant bracket does not satisfy the Jacobi identity in spite of the fact that there is a striking similarity between this bracket and the Lie bracket on the double ${\cal G} \oplus {\cal G}^*$ of a 
 Lie bialgebra  $({\cal G}, {\cal G}^*)$ (see [Dr]). In fact, if we denote by
 $ \lag .,. \rg_+$ the canonical symmetric 2-form
 on $TM \oplus T^*M$, this induces an operation on $\Gamma(TM \oplus T^*M)$ 
  defined by 

$$ \lag X_1+ \alpha_1, X_2+ \alpha_2 \rg_+ ={1 \over 2}(i_{X_2} \alpha_1 + i_{X_1} \alpha_2).$$

\noindent Then,  

$$[[e_1,e_2]_c, e_3]_c \ +c.p.={1 \over 3}d \lag [e_1,e_2]_c, e_3 \rg_+ \ +c.p. $$

\noindent for any $e_1,$ $e_2$,  $e_3 \in \Gamma(TM \oplus T^*M)$. 
 Here the notation $c.p.$ denotes the other terms obtained by cyclic permutations  of indices  1,2 and 3. The Courant bracket was systematized  by Liu,
  Weinstein and Xu in~\cite{LWX}. \par
  In what follows, for simplicity, we will
   omit the base space of  vector bundles
 when this is implied (i.e. a vector bundle $E \rightarrow M$ 
 will be often replaced by its total space $E$).

\begin{dfn}{\rm
Let $L$ be a sub-bundle of $TM \oplus T^*M$. Then $L$ is said to be a 
 {\sl Dirac structure on the manifold} $M$ if $L$ is maximally isotropic 
 for the bilinear form $\lag .,. \rg_+$ and the set $\Gamma(L)$ of 
 smooth sections of $L$  is closed under the Courant bracket.
}
\end{dfn}

 \noindent Let us give the basic examples of Dirac structures:

\begin{ex} {\bf Pre-symplectic 2-forms.} \par
{\rm  Let $ \omega$ be a differential 2-form on $M$. It can be considered 
as a map $\omega: TM \rightarrow  T^*M$ defined by 
$\omega(Y)=i_Y \omega$. Let us denote by $L_{\omega}$ its graph.
Then, $L_{\omega}$ is a Dirac structure on $M$ if and only if $\omega$ is a closed 2-form.
}
\end{ex}

\begin{ex} {\bf Poisson structures.}\par
{\rm  A {\sl Poisson structure} on a manifold $M$ is given by 
 a bivector field $\pi \in \Gamma(\Lambda^2(TM))$  such that $[\pi,\pi]_s=0$,  where $[ \ , \ ]_s$ is the Schouten-Nijenhuis bracket on the space of 
 poly-vector fields (see~\cite{Kz}). \par

 Any bivector field $\pi$ defines a skew symmetric map, that we denote also 
 by $ \pi: T^*M \rightarrow  TM$, whose extension to the sections is 
  given by $i_{\pi (\alpha)} \beta= \pi(\alpha, \beta)$ for any two 1-forms 
 $\alpha$, $\beta$.  The graph $L_{\pi}$ of that map is a Dirac structure on 
 $M$ if and only if $\pi$ is a Poisson structure.
}
\end{ex}

\begin{ex} {\bf Involutive distributions.}\par

{\rm Let $F$ be sub-bundle of $TM$ and denote by $F^{\perp}$ its annihilator in $T^*M$. We consider $L=\{Y+\xi \ | \ Y \in F, \quad \xi \in F^{\perp}\}.$ 
  Then $L$ is a Dirac structure on $M$
if and only if $\Gamma(F)$ is closed under the Lie bracket of vector fields.
}
\end{ex}

\section{ Extension of the Courant Bracket }

 Let $M$ be a smooth manifold. We shall extend the Courant bracket to the smooth sections of the
  vector bundle ${\cal E}^1(M)\rightarrow M$, where

$${\cal E}^1(M)= (TM \times \reals) \oplus (T^*M \times \reals).$$

\noindent The smooth sections of ${\cal E}^1(M)$ are of the form 
 $e=(X,f)+(\alpha,g)$, where 
$X$ is a smooth vector field on $M$, $\alpha$ is a  1-form and $f$, $g$  are
smooth functions on $M$. In what follows, we use the notation  $X \cdot f$ 
 instead of $i_Xdf$ for any vector field $X$ and for any smooth function $f$.
  Consider the following bracket on $\Gamma({\cal E}^1(M)$):
\begin{eqnarray*}
 \lefteqn [ \ (X_1, f_1)+(\alpha_1, g_1),&(X_2, f_2)+(\alpha_2, g_2) ] = \Big([X_1,X_2], \  X_1  \cdot f_2 - X_2 \cdot f_1 \Big) \ \ \ \ \ \  \ \ \ \ \ \ \ \ \ \ \cr
 & +\Big( \call_{X_1} \alpha_2-\call_{X_2} \alpha_1
 +{1 \over 2}d(i_{X_2}\alpha_1- i_{X_1}\alpha_2) \cr
& + f_1 \alpha_2 -f_2 \alpha_1 +{1 \over 2}(g_2df_1-g_1df_2 - f_1dg_2+f_2dg_1 ),  \cr
&  X_1 \cdot g_2- X_2 \cdot g_1  
 +{1 \over 2}(i_{X_2}\alpha_1 - i_{X_1}\alpha_2 -f_2g_1+f_1g_2) \Big),
\end{eqnarray*}

\noindent where  the $(X_i,f_i)+(\alpha_i,g_i)$ are in $\Gamma({\cal E}^1(M))$.
 In order to simplify our notation, we use the same symbol 
 $[.,.]$ for both the Lie bracket on vector fields and the bracket on 
 sections of ${\cal E}^1(M)$. By identifying $(X_i,0)+(\alpha_i,0)$
 with $X_i+\alpha_i$, we get: 

\begin{eqnarray*}
[X_1+\alpha_1, X_2+\alpha_2]& = \Big([X_1,X_2],0 \Big)+\Big(\call_{X_1} \alpha_2 -\call_{X_2} \alpha_1+{1\over 2}d(i_{X_2}\alpha_1- i_{X_1}\alpha_2), \cr
 & \hfill {1 \over 2}(i_{X_2}\alpha_1- i_{X_1}\alpha_2)\Big).
\end{eqnarray*}

\noindent So the projection of the second member of this equality onto $TM \oplus T^*M$ gives the Courant bracket of $X_1+\alpha_1$ and $X_2+\alpha_2$.

\begin{rmk}
\label{injection}
{\rm To any sub-bundle $L \subset TM \oplus T^*M$, we can associate 
the sub-bundle $\widetilde L$ of ${\cal E}^1(M)$ whose fibre at $x$ is given by 
$$\widetilde L_x = \{ (X(x),0)+(\alpha(x), f(x))  \ | \  X+\alpha \in \Gamma( L),
 \ f \in C^{\infty}(M)\}.$$

\noindent In fact $L$ is maximally isotropic for the 2-form $ \lag .,. \rg_+$ 
 if and only if $\widetilde L$ is a maximal isotropic sub-bundle
 with respect to the extension of $ \lag .,. \rg_+$ to ${\cal E}^1(M)$ 
 defined by:

$$\Big \lag(X_1, \mu_1)+ (\alpha_1, \lambda_1), (X_2, \mu_2)+ (\alpha_2, 
\lambda_2)
 \Big \rg_+ 
={1 \over 2}\Big(i_{X_2} \alpha_1 + i_{X_1} \alpha_2
+\lambda_1 \mu_2 + \lambda_2 \mu_1\Big)$$

\noindent for any $(X_i,\mu_i)+ (\alpha_i,\lambda_i)$ elements of 
 ${\cal E}^1(M)$, with $i=1,2$. Furthermore, $\Gamma(L)$ is closed under the Courant bracket if
and only if  $\Gamma(\widetilde L)$ is closed under the bracket $[.,.]$.

} \end{rmk}

\begin{dfn}{\rm
Let $L$ be a sub-bundle of ${\cal E}^1(M)$. Then $L$ is be said to be a 
{\sl Dirac structure of} ${\cal E}^1(M)$ if $L$ is maximally isotropic for 
 the bilinear form $ \lag .,. \rg_+$ and $\Gamma(L)$  is closed under 
  $[.,.]$.
}
\end{dfn}

 It follows from Remark~\ref{injection} that Dirac structures on $M$ 
  are in one-to-one correspondence with Dirac structures of ${\cal E}^1(M)$
  which are subsets of $TM \oplus (T^*M \times \reals)$. 
 Hence, any closed 2-form, foliation or  Poisson structure
  can be identified with a Dirac structure of ${\cal E}^1(M)$.

Let us fix some notation.

\smallskip
\noindent{ \bf Notations.} For  any $e_1,e_2,e_3 \in \Gamma({\cal E}^1(M)),$
 we set

$$T(e_1,e_2,e_3)= {1 \over 3} \lag[e_1,e_2],e_3 \rg_+ +c.p
 \quad {\rm and} \quad \rho(e)h= X \cdot h$$

\noindent  for any $e=(X, f)+(\alpha,g) \in \Gamma({\cal E}^1(M))$ and 
 $h \in C^{\infty}(M)$. 

\begin{prop} 
\label{properties}
{\sl With the above notations, we have:

\begin{itemize}
\item[(i)] \ $[[e_1,e_2],e_3]+c.p.=(0,0)+\Big(dT(e_1,e_2,e_3), \ T(e_1,e_2,e_3)\Big)$ for any $e_1,e_2 ,e_3 \in \Gamma({\cal E}^1(M)), $
\item[(ii)] \  $[e_1,fe_2]=f[e_1,e_2]+ \big(\rho(e_1)f \big)e_2 
 -\lag e_1,e_2 \rg_+ \Big((0,0)+(df,0) \Big)$  for any smooth function $f$ and for any  $e_1, e_2 \in \Gamma({\cal E}^1(M))$.
\end{itemize}
}
\end{prop}

\noindent The proof of Proposition~\ref{properties} is straightforward.
 We shall use this proposition to show that any Dirac structure of 
 ${\cal E}^1(M)$ is a  Lie algebroid. Recall that a vector bundle
  ${\cal E}$ over a smooth manifold $M$ is said to be a 
 {\sl Lie algebroid} if there is a Lie bracket $[ \ , \ ]_{\cal E}$ on 
  $\Gamma({\cal E})$ and a  bundle map 
 $\varrho : {\cal E}  \rightarrow TM$, extended to a map 
 between sections of these bundles, such that

\begin{itemize}
\item[1)] $\varrho([X, Y]_{\cal E})=[\varrho(X),\varrho(Y)],$
 
\item[2)] $[X, fY]_{\cal E}= f[X,Y]_ {\cal E}+ (\varrho(X)f)Y$,
\end{itemize}

\noindent for any $X$, $Y$ smooth sections of ${\cal E}$ and for any smooth
  function $f$  defined on $M$. 
 Then $\varrho$ is called the {\sl anchor} of the Lie algebroid.

\begin{thm}
\label{algebroid}
 {\sl Let $L$ be an isotropic sub-bundle of $({\cal E}^1(M), \lag .,. \rg_+)$
  such that $\Gamma(L)$ is closed under the bracket  $[.,.]$ of
 $\Gamma({\cal E}^1(M))$. Then $(L, [.,.], \rho_{|_L})$ is a Lie algebroid.
 In particular, any Dirac structure of ${\cal E}^1(M)$ is a Lie algebroid.
 }
\end{thm}

\noindent{\sl Proof:} By applying the first property of Proposition~\ref{properties},  we obtain the Jacobi identity on $\Gamma(L)$. Moreover, the second 
property of Proposition~\ref{properties} implies that 
$$[e_1,fe_2]=f[e_1,e_2]+ (\rho(e_1)f)e_2 \quad \forall e_1, e_2 \in \Gamma(L),
\quad  \forall f \in C^{\infty}(M).$$
\hfill \qed

 \begin{rmk}{\rm
Let $(A, A^*)$ be a Lie bialgebroid (i.e. $A$ and $A^*$ are endowed with Lie
 algebroid structures such that the differential $d_*$ on $\Gamma(\wedge A)$,
 induced from the Lie algebroid structure on A$^*$ is a derivation of the
 Lie bracket on $\Gamma(A)$). In~\cite{LWX}, the authors define a bracket
 on the direct sum ${\cal E}=A\oplus A^*$ such that ${\cal E}$ is a 
 {\sl Courant algebroid} (see~\cite{LWX} for details). This bracket is different 
 from the extended Courant bracket that we define above. It is clear that 
 $(TM, T^*M)$ is a bialgebroid, where $\Gamma(TM)$ is equipped with the Lie 
 bracket and  $T^*M$ is equipped with the null Lie algebroid bracket.
  But Property (i) of Proposition~\ref{properties} shows that the 
 extended Courant bracket defined on $\Gamma({\cal E}^1(M))$ does not induce
  a Courant algebroid structure on ${\cal E}^1(M)$.   
}
\end{rmk}

\section { Further Examples of Dirac Structures}

\noindent We begin this section by giving a correspondence
 between Jacobi manifolds and some special Dirac structures.
  Next, we will show that apart from 
 Jacobi manifolds, there are other interesting examples of Dirac structures.

\subsection {Jacobi Manifolds}
 A {\sl Jacobi structure}  on a manifold $M$
  is given by a pair $(\pi, E)$ formed by a bivector field $\pi$
   and a vector field $E$ such that (~\cite{L})

$$[E,\pi]_s=0 \quad {\rm and} \quad [\pi, \pi]_s=2E \wedge \pi,$$ 

\noindent where $[ \ , \ ]_s$ is the Schouten-Nijenhuis bracket on the space of
poly-vector fields. A manifold endowed with a Jacobi structure is said
 to be a {\sl Jacobi manifold}.  When $E$ is zero, we get a Poisson structure. 
  \par

 Any smooth section $P$ of $\G^2(TM \times R)$ is defined by a pair $(\pi, E)$, where $\pi$ is a bivector field and $E$ is a vector field on $M$. Indeed,
   $P$  can be regarded as  a vector bundle map 
 $ P$: $T^*M \times R \rightarrow TM \times R$  whose extension to the sections of $T^*M \times R $ has the following form:
 $$ P(\alpha, g)= (\pi \alpha+ g E, -i_E\alpha),$$

\noindent for any 1-form $\alpha$ and for any smooth function $g$, where
 $\pi \alpha$ is the vector field defined by $i_{\pi \alpha} \beta= \pi(\alpha, \beta)$ for any 1-form $\beta$.  The graph  $L_{ P}$ of $P$ is a Dirac structure of ${\cal E}^1(M)$ 
 if and only if $(\pi,E)$ is a Jacobi structure. In such a case, we can
 identify $L_{ P}$ with $T^*M \times \reals$ equipped with the Lie algebroid 
 structure that is associated to the Jacobi manifold.
 Namely, we consider the map 
  $\Phi: T^*M \times \reals \rightarrow L_{ P}$ defined by 
 $\Phi(\alpha, f)= (\pi \alpha+ g E, -i_E\alpha)+ (\alpha, f)$.
Then $\Phi$ a Lie algebroid isomorphism.
  Here  $\Gamma(L_{ P})$ is endowed with the extended bracket $[.,.]$, while
 $\Gamma(T^*M \times \reals)$ is equipped with the
following bracket introduced in~\cite{KS}

\begin{eqnarray*}
\{(\alpha,f),(\beta ,g)\}_{_{(\pi, E)}}&=\Big({\cal L}_{\pi \alpha}\beta -
 {\cal L}_{\pi \beta}\alpha -d(\pi(\alpha, \beta)) +f {\cal L}_E\beta -g
 {\cal L}_E\alpha -i_E(\alpha \wedge \beta), \\
&  -\pi(\alpha, \beta) + \pi(\alpha, dg)-\pi(\beta, df) +fE \cdot g -gE \cdot f \Big).
\end{eqnarray*}

\noindent We refer the reader to [H-M] for discussions on morphisms of 
 Lie algebroids. One can observe that  $(T^*M \times R, TM \times R)$
 is not a Lie bialgebroid when $(M,\pi,E)$ is a generic Jacobi manifold
 (otherwise we obtain the Poisson case).

\subsection {Nambu Manifolds.}

 Let $M$ be an $n$-dimensional smooth manifold.
   A {\sl Nambu structure} on $M$ of order $p$ ($2\leq p \leq n$) 
 is defined by a $p$-vector field  $\Pi$ which satisfies the Fundamental 
 Identity (see~\cite{T})

\begin{eqnarray*}
 \{f_1,...,f_{p-1}, \{g_1,...,g_p\}_{_\Pi} \}_{_\Pi} =
 \sum_{k=1}^{p}  \{g_1,...,g_{k-1}, \{f_1,...,f_{p-1},g_k\}_{_\Pi},
 g_{k+1},..., g_p\}_{_\Pi}, 
\end{eqnarray*}

\noindent  for any $ f_1,\dots ,f_{p-1}, g_1, \dots g_p \in C^{\infty}(M)$,
  where $\{ \ \ \}_{_\Pi}$ is  defined by 

$$\{f_1, \dots ,f_p \}_{_\Pi}=\Pi(df_1, \dots , df_p), \quad \forall
 \ f_1,\dots ,f_{p} \in C^{\infty}(M).$$ 

 \noindent A manifold equipped with such a structure is called 
 {\sl Nambu manifold}. In fact, Nambu structures of order 2 are nothing 
 but Poisson structures.  We shall show that the local structure of
 Nambu manifolds gives rise to certain Dirac structures:

\begin{prop} \label{nambu-dirac}
{\sl
  Locally, any Nambu structure of order $n-2$ on an $n$-dimensional
  manifold corresponds to a family of Dirac structures.
 }
 \end{prop}

 To prove this proposition, we shall use the characterization of Nambu 
 manifolds by means of differential forms which is given in~\cite{DZ}. 
 Precisely, let $\Omega$ be
 a volume form on $M$,  any $p$-vector field $\Pi$ corresponds to
 a $(n-p)$-form $\omega= i_{\Pi} \Omega$. In~\cite{DZ}, the authors
  proved that $\Pi$ is a Nambu structure if and only if 
 
\begin{equation}
\label{eq:co-Nambu}
 (i_A \omega) \wedge \omega =0 \quad {\rm and} \quad
  (i_A \omega) \wedge d\omega =0,
\end{equation}

\noindent for any $(n-p-1)$-vector field $A$. 
 Obviously, $\omega$ depends on the
 volume form $\Omega$. A differential form which satisfies
  the relations (~\ref{eq:co-Nambu}) is called a {\sl co-Nambu form}.
 For the proof of  Proposition~\ref{nambu-dirac},
 we need  also the following two lemmas:

 \begin{lemma}(see~\cite{Gaut})
 \label{gauth}
{\sl A $p$-vector field $\Pi$ is a Nambu structure of order $p$ if and only if
 for any point $x$ where $\Pi(x) \ne 0$, there exists a local system 
 of coordinates $(x_1, \dots, x_n)$ defined in a neighborhood of $x$ such that 
  $$\Pi= {\partial \over \partial x_1} \wedge \dots \wedge
  {\partial \over \partial x_p}.$$
 }
 \end{lemma}

\begin {lemma}\label{struct}
 {\sl Any pair $(\omega, \eta)$ formed by a differential 2-form and
 a closed 1-form such that $d \omega= \eta  \wedge \omega$
  corresponds to a Dirac structure $L_{(\omega, \eta)}$ of ${\cal E}^1(M)$.
}
 \end {lemma}
\noindent {\sl Proof:}
 For any  2-form $\omega$ and any closed 1-form $\eta$, we
 consider the sub-bundle $L_{(\omega, \eta)} 
 \subset {\cal E}^1(M)$  whose fibre at a
 point $x$ is:
 
$$L_{(\omega, \eta)}(x)=\{ (X, -i_X \eta)_x+(i_X \omega + f \eta, \ f)_x  \ |  \ X \in \Gamma(TM), \ f \in C^{\infty}(M) \}.$$

\noindent It is not hard to check that $L_{(\omega, \eta)}$ is a 
 maximal isotropic sub-bundle of ${\cal E}^1(M)$.  Furthermore, 
 by a simple (but long) computation, we prove that  
 $L_{(\omega, \eta)}$ is a Dirac structure of ${\cal E}^1(M)$ if and only if
  we have the relation 

 \begin{equation}
\label{eq:strano}
 d \omega= \eta  \wedge \omega.
\end{equation}

\noindent This gives the lemma.
\hfill \qed

\bigskip
\begin{rmk} {\rm Dirac structures of the type $L_{(\omega, \eta)}$
  include {\sl locally conformal symplectic structures}. Recall that   
 a locally conformal symplectic structure on an $2m$-dimensional
 manifold $M$ is given by a pair $(\omega, \eta)$, where $\omega$
 is a non-degenerate 2-form and $\eta$ is a closed 1-form such that
 $d \omega= \eta  \wedge \omega$. We have the definition:
}
\end{rmk}

\begin{dfn}{\rm A {\sl locally conformal pre-symplectic structure} is given
by a 2-form $\omega$ and a closed 1-form such that 
$d \omega= \eta  \wedge \omega$.
}
\end{dfn}

 \noindent {\sl Proof of Proposition~\ref{nambu-dirac}:}  
 Let $\Pi$ be a Nambu structure of order $n-2$ on an $n$-dimensional manifold.
  Denote ${\cal U}= \{ x \in M \ | \ \Pi(x) \ne 0\}$. 
 Taking into account the fact that $\Pi$ is a Nambu structure
 on $M$ if and only if it is such on ${\cal U}$, we can restrict ourselves
 to  ${\cal U}$. According to Lemma~\ref{gauth}, around any point $x_0$ 
 of ${\cal U}$ we can write: 
  $\Pi= \partial / \partial x_1 \wedge \dots \wedge
  \partial /  \partial x_{n-2}$  for a suitable system of coordinates
   $(x_1, \dots, x_n)$ defined on an open set ${\cal V}_{x_0}$ containing
 $x_0$.  To any function $f$ such that $f(x) \ne 0$ on ${\cal V}_{x_0}$,
 we associate the local volume form
  $\Omega_f= f dx_{1} \wedge \dots \wedge dx_n$. We
 obtain  $i_{\Pi} \Omega_f = \pm f dx_{n-1} \wedge dx_n$. Denote 
 $\omega_f=i_{\Pi} \Omega_f$. Then we have
   $d \omega_f = \pm dln|f| \wedge \omega_f$.
  Applying Lemma~\ref{struct}, 
 we deduce the existence of a Dirac structure of $E^1({\cal V}_{x_0})$
  associated to any pair $(\omega_f, dln|f|)$. Thus, any
 function $f: {\cal V}_{x_0} \rightarrow \reals$ that vanishes 
 nowhere,  defines a Dirac structure $L_{(\omega_f, dln|f|)}$. \par

  Conversely, assume that $L_{(\omega, dln|f|)}$ is Dirac structure,
  where $f$ is a function which does not vanish on $M$ 
  and $\omega$ is a decomposable 
  2-form (i.e.  there are two independent differential 1-forms 
 $\omega_1$, $\omega_2$ such that $\omega=\omega_1 \wedge \omega_2$). 
 Then, $\omega$ is a co-Nambu 2-form. Moreover,  $g \omega$ is also a co-Nambu
  2-form for any non-vanishing function $g$. From
  Lemma~\ref{struct}, we obtain that any pair $(g \omega, dln|fg|)$ 
 corresponds to a Dirac structure.
  Using a local volume form $\Omega$ defined on an open set ${\cal V}$, 
 we get a Nambu structure $\Pi$ on ${\cal V}$ characterized by 
  $i_{\Pi} \Omega= \omega$.
 Therefore,  Proposition~\ref{nambu-dirac} is obtained.

 \hfill \qed

\subsection { Homogeneous Poisson Manifolds.}

A {\sl homogeneous Poisson manifold} $(M,\pi, Z)$ is a Poisson manifold
  $(M,\pi)$ with a vector field $Z$ on $M$ such that $[Z,\pi]_s=-\pi$,
 where $[.,.]_s$ is the Schouten-Nijenhuis bracket defined on poly-vector 
 fields. This kind of structures was studied by Dazord,
  Lichnerowicz and Marle in~\cite{DLM}. We shall give a new characterization of these structures. Namely, we have the following result:

\begin{prop}
\label{ex:homo}
 {\sl Let $\pi$ be a bivector field and let $Z$ be a vector field on $M$. 
Consider the sub-bundle $L_{(\pi, Z)}$ whose fibre at a point $x$ is given
$$L_{(\pi, Z)}(x)=\{ (\pi \alpha -fZ, \ f)_x+(\alpha, \ i_Z \alpha)_x  
 \ |  \ \alpha \in \Gamma(T^*M), \ f \in C^{\infty}(M) \}.$$

\noindent Then $(M,\pi, Z)$ is a homogeneous Poisson manifold if and only if
 $L_{(\pi, Z)}$ is a Dirac structure of ${\cal E}^1(M)$. 
 }
\end{prop}

\noindent {\sl Proof:} By an easy computation, we see
  that $L_{(\pi, Z)}$ is maximally isotropic. 
 Notice that the set of smooth sections of  
 $L_{(\pi, Z)}$ is spanned by $(Z, -1)+(0,0)$ and elements of the type  
 $(\pi \alpha,0)+(\alpha, \ i_Z \alpha)$, where $\alpha$ is an arbitrary 
  differential 1-form. By using Property $(ii)$ of 
 Proposition~\ref{properties}, we see that 
 it is sufficient to study the bracket of such sections.
  On one hand for any differential 1-form,  we have
 
\begin{eqnarray*}
 [(\pi \alpha,0)+(\alpha, \ i_Z \alpha), \ (Z, -1)+(0,0)]&= &-
 \Big([Z, \pi \alpha] + \pi {\cal L}_Z \alpha,\  0\Big)\cr
& {} & -\Big({\cal L}_Z \alpha - \alpha, \ i_Z({\cal L}_Z \alpha- \alpha)\Big).
\end{eqnarray*}
 
\noindent Using the definition of the Schouten-Nijenhuis bracket, we get 

\begin {equation}\label{eq:formula}
[Z, \pi]_s \ \alpha=[Z,\pi \alpha]-\pi {\cal L}_Z \alpha \quad \forall \alpha \in 
\Gamma(T^*M). \end{equation} 
\noindent It follows that $[(\pi \alpha,0)+(\alpha, \ i_Z \alpha),
 \ (Z, -1)+(0,0)]$  belongs to $\Gamma(L_{(\pi, Z)})$
 if and only if 
 $$[Z, \pi]_s \ \alpha =-\pi \alpha.$$

\noindent On the other hand, from a direct calculation of the bracket

$$[(\pi \alpha_1,0)+(\alpha_1, \ i_Z \alpha_1), \ 
 (\pi \alpha_2,0)+(\alpha_2, \ i_Z \alpha_2)],$$
\noindent we deduce that this bracket is a smooth section of 
  $L_{(\pi, Z)}$ for any two differential 1-forms $\alpha_1$, $\alpha_2$ 
 if and only if $[\pi, \pi]_s=0$ and $[Z, \pi]_s=-\pi$.

\noindent This completes the proof of Proposition~\ref{ex:homo}.

\hfill \qed

\smallskip

 Before ending this section, let us remark that any exact differential 2-form 
 on a manifold $M$ is a Dirac structure of ${\cal E}^1(M)$ which
  is a graph of a certain map. More precisely, let $\mathbold B$ be a smooth section 
 of $\G^2(T^*M \times R)$. So,  $\mathbold B$ is given by a pair 
 $(\omega, \alpha)$ formed by a 2-form $\omega$ and a 1-form $\alpha$.
  We can define a vector bundle map
$\phi_{\mathbold B} :TM \times R \rightarrow T^*M \times R$ extended to 
sections, such that

$$\phi_{\mathbold B}(X,f)= (i_X \omega+f \alpha, -i_X \alpha) $$

\noindent for any smooth vector field $X$ and for any $f \in C^{\infty}(M)$.
 The graph of $\phi_{\mathbold B} $ is a Dirac structure of ${\cal E}^1(M)$ 
  if and only if $ \omega=d \alpha$.

\section {The Jacobi Algebra of Admissible Functions}

\begin{dfn}{\rm
Let $L$ be a maximal isotropic sub-bundle of $({\cal E}^1(M), . \lag .,. \rg)$. A function $f$ is said to be $L$-{\sl admissible} if there exist a vector field $X_f$ and a smooth function 
 $\varphi_{f}$ on $M$ such that  $e_f=(X_f,\varphi_{f})+(df,f)$ is a
  smooth section of $L$. 
}
\end{dfn}

   Note that $e_f$ is unique up to a smooth section of 
 $L \cap (TM \times \reals)$. Let $ \lag .,. \rg_-$ denote the skew symmetric
  2-form on ${\cal E}^1(M)$ defined by:

$$\Big \lag(X_1, \mu_1)+(\alpha_1, \lambda_1), (X_2, \mu_2)+(\alpha_2, 
\lambda_2)  \Big \rg_-  = {1 \over 2} \Big(i_{X_2} \alpha_1 - i_{X_1}\alpha_2
+\lambda_1 \mu_2 -\lambda_2 \mu_1\Big)$$

\noindent for any $(X_i, \mu_i)+(\alpha_i, \lambda_i)$ in  ${\cal E}^1(M)$.
 Define the following bracket on the space of  $L$-admissible functions:

$$\{f,g\}= - \lag e_f,e_g \rg_-$$

\noindent  for any two $L$-admissible functions $f$ and $g$. This bracket is well defined. Indeed, if 
$e'_f=(X'_f,\varphi'_{f})+(df,f)$  is another element of $\Gamma(L)$ then 
$e_f-e'_f$ is in  $\Gamma(L\cap (TM \times \reals))$. 
 Therefore for any $L$-admissible function $g$, we have
$$\Big \lag e_f-e'_f, \ e_g \Big \rg_+ =
 \Big \lag  e_f-e'_f,\  e_g\Big \rg_-=0.$$ 

\noindent Hence $\lag e_f,e_g \rg_-=\lag e'_f,e_g \rg_-.$  Since
$\lag e_f,e_g \rg+=0$, we can rewrite 
 $$\{f,g\}= i_{X_f}dg+g \varphi_f= - (i_{X_g}df+f \varphi_g).$$ 
 It should be mentioned that the bracket $\{.,.\}$ is local, this means
that the support of $\{f,g\}$ is contained in the intersection of the supports of $f$ and $g$. Precisely,  $\{f,g\}(x)$ depends only on the 
values of $df(x)$, $dg(x)$, $f(x)$ and $g(x)$.
 One can notice that our definition of an $L$-admissible function is a 
 little bit more general than the one given by Courant in~\cite{C}.
 Here, we do not need  the fact that $\Gamma(L)$ is closed under $[.,.]$ in
 the definition of $\{.,.\}$. We have the following result:

\begin{prop}
\label{algebra}
 {\sl Let $L$ be a Dirac structure of ${\cal E}^1(M)$. Then the space of
  $L$-admissible functions is a Lie algebra.
 }
 \end{prop}

\begin{lemma}
\label{prep1} 
 {\sl Let $L$ be an isotropic sub-bundle of ${\cal E}^1(M)$ for the 
 symmetric form 
 $ \lag .,, \rg_+.$ If  $e_f=(X_f,\varphi_{f})+(df,f)$ and 
 $e_g=(X_g,\varphi_{g})+(dg,g)$ are two smooth sections of $L$, then 
$$[e_f,e_g]=([X_f,X_g], \ X_f \cdot \varphi_{g}- X_g \cdot \varphi_{f})
+(d\{f,g\},\{f,g\}),$$

\noindent where $\{f,g\}= -\lag e_f,e_g \rg_-.$
 }
 \end{lemma}

\noindent {\sl Proof:}  Using the definition of the bracket $[.,.]$ on
 ${\cal E}^1(M)$, we get 
\begin{eqnarray*}
 [e_f,e_g]&=([X_f,X_g], \ X_f \cdot \varphi_{g}- X_g \cdot \varphi_{f}) 
 +{1 \over 2}\Big (d(X_f \cdot g - X_g \cdot f+g\varphi_{f} - f\varphi_{g}),\cr
 & {}  \quad \quad \quad \quad \quad  \quad  \quad  \quad \quad \quad \quad
  X_f \cdot g - X_g \cdot f+g\varphi_{f} - f\varphi_{g} \Big).
\end{eqnarray*}
\noindent  This is equivalent to 
\begin{eqnarray*}
[e_f,e_g]&=& ([X_f,X_g], \ X_f \cdot \varphi_{g}- X_g \cdot \varphi_{f})-
 (d\lag e_f,e_g \rg_-, \lag e_f,e_g \rg_-)\cr
&= & ([X_f,X_g], \ X_f  \cdot\varphi_{g}- X_g \cdot \varphi_{f})
+(d\{f,g\},\{f,g\}).
\end{eqnarray*}
\noindent This proves the lemma.

\smallskip

\hfill \qed

\noindent {\sl Proof of Proposition~\ref{algebra}:}
 Let $L \subset{\cal E}^1(M)$ be a Dirac structure. It is clear that the 
 corresponding bracket $\{.,.\}$ is $\reals$-bilinear and skew symmetric. Moreover, using Lemma~\ref{prep1},
 we get:

 \begin{eqnarray*}
 \lefteqn \lag \ [e_f,e_g], \ e_h\rg_+&= {1 \over 2} \Big(i_{[X_f,X_g]} dh + h(\ X_f \cdot \varphi_{g}- X_g \cdot \varphi_{f}) + \{h, \{f,g\} \} \Big)  \ \ \ \ \ \ \ \ \ \ \ \ \ \ \ \ \ \ \ \ \cr
&={1 \over 2} \Big({\cal L}_{X_f}{\cal L}_{X_g}h -
{\cal L}_{X_g}{\cal L}_{X_f}h + h(\ X_f \cdot \varphi_{g}- X_g \cdot \varphi_{f}) + \{h, \{f,g\} \} \Big)
\end{eqnarray*}
 
\noindent for any $L$-admissible functions $f$, $g$ and $h$. 
If we add and withdraw the term $\varphi_{f} X_g \cdot h -
\varphi_{f} X_g \cdot h$ in the second member of this last equality, we obtain
\begin{eqnarray*}
\lag \ [e_f,e_g], \ e_h\rg_+&={1 \over 2} \Big(X_f \cdot \{g,h\} 
 -X_g \cdot \{f,h\}  + \varphi_{f} X_g \cdot h -\varphi_{g} X_f \cdot h + \{h,\{f,g\}\} \Big)\cr
 &=\{f, \{g,h\} \}-\{ g, \{f,h\} \}+ \{h,\{f,g\} \}.
\end{eqnarray*}

\noindent  By using the definition of a Dirac structure, we have:

$$\Big \lag [e_f,e_g], \ e_h \Big \rg_+=\{f, \{g,h\} \}-\{ g, \{f,h\} \}+
\{h,\{f,g\} \}=0.$$

\noindent Therefore, the space of $L$-admissible functions forms a Lie algebra.

\hfill \qed

\smallskip
 
 Proposition~\ref{algebra} extends a result proven in~\cite{CW} for Dirac structures on $M$ (i.e. for Dirac structures of ${\cal E}^1(M)$  which  are
 subsets of $TM \oplus (TM \times \reals$)). Now let us give a remark that will
 be useful later.

\begin{rmk}
\label{product}
{ \rm Let $L$ be a Dirac structure of ${\cal E}^1(M)$. 
 The product $fg$ of two $L$-admissible functions $f$ and $g$ may not be 
 $L$-admissible. But, if
 $(Y_g, \theta_g)+(dg,0)$ is in $\Gamma(L)$, then for any $L$-admissible 
 function $f$, the product $gf$ is also $L$-admissible. Indeed, the section

$$ g \Big( (X_f, \varphi_f)+(df,f) \Big)+f \Big((Y_g, \theta_g)+(dg,0) \Big)=
(gX_f+ f Y_g, \ g \varphi_{f}+f \theta_g)+( d(gf),gf)$$

\noindent belongs to $\Gamma(L)$. By induction,  we obtain that
$g^n f$ is $L$-admissible for any integer $n \geq 1$. 
On the other hand, If $f$, $g$, $h$ and $gh$ are $L$-admissible functions, then we have 

\begin{equation}
\label{eq:leibniz defaut}
 \{f,gh\}=g\{f,h\} +\{f,g\}h -gh \ \varphi_f,
\end{equation}

\noindent where $(X_f, \varphi_f)+(df,f) \in \Gamma(L)$. This holds since 
 we have:
\begin{eqnarray*}
\{f,gh\} &=&g X_f \cdot h + h X_f \cdot g +gh \varphi_f \cr
&=& g (X_f \cdot h+ h \varphi_f)+(X_f \cdot g+ g \varphi_f)h -gh \varphi_f\cr
&=&g\{f,h\} +\{f,g\}h -gh \ \varphi_f.
\end{eqnarray*}

\hfill \qed

}
\end{rmk}

 From now on, we denote by $A^{\infty}(L)$ the set of $L$-admissible functions,
 whenever $L$ is maximal isotropic sub-bundle of 
 $({\cal E}^1(M), \lag .,. \rg_+)$.\\

 The rest of this section is devoted to the case where the constant functions
  are $L$-admissible, for a given Dirac structure $L \subset {\cal E}^1(M)$.
  In such a case, there exists a smooth vector field $E$ such that 
 $(E,0)+(0,1)$ belongs to $\Gamma(L)$. Therefore, if $e_f=(X_f, \varphi_f)+(df,f)$ is in $\Gamma(L)$, then
 $$\lag (E,0)+(0,1), \ e_f \rg_+=0.$$

\noindent This gives  $\varphi_f= -E \cdot f$. Hence 
 Equation (~\ref{eq:leibniz defaut}) implies

$$\{f,gh\}=g\{f,h\} +\{f,g\}h +gh \ E \cdot f.$$
It follows from this above equality, the existence a bi-differential operator
  on $A^{\infty}(L)$ of order
 at most one in each of its argument, denoted by $P$, such that
  $P(f,g)=\{f,g\}$.
 This shows that $A^{\infty}(L)$ is a Jacobi algebra. \par

Recall that a {\sl Jacobi algebra} 
is an associative  commutative algebra ${\cal A}$ with unit 1 (over $\reals$ 
 or ${\mathbb C}$) endowed with a skew symmetric bi-differential operator $P$, 
 which defines a Lie algebra structure on ${\cal A}$. When $P$ vanishes on constants, we get a {\sl Poisson algebra}.
 A Jacobi structure $(\pi, E)$ on a smooth manifold $M$ is equivalent to 
 having a Jacobi  algebra structure on $C^{\infty}(M)$ (~\cite{Gr}). 
 We summarize our discussion in the following form:

\begin{thm}
 {\sl There exists an one-to-one correspondence between Jacobi structures 
 on $M$ and Dirac structures $L$ of ${\cal E}^1(M)$ such that 
 $A^{\infty}(L)=C^{\infty}(M)$.
}
\end{thm}

\section {Conformal Dirac Structures}

In this section, we introduce the concept of a conformal Dirac structure
  which includes conformal Jacobi structures. Conformal Jacobi structures 
 are known to be natural generalizations of local Lie algebras considered by 
 Kirillov in~\cite{K}. We refer the reader to~\cite{DLM}, for instance, for 
  details on conformal Jacobi structures.
  First, we establish the following result:

\begin{thm}
\label{step2} {\sl Let $L$ be a Dirac structure of ${\cal E}^1(M)$. Assume that
  $\Gamma(L)$ contains an element of the form $(Y_a, \theta_a)+(da,0)$, 
 with $a(x) \ne 0$ for any  $x \in M$.  Then the space 
 $A^{\infty}(L)$  of $L$-admissible functions is endowed with two
  Lie algebra brackets $\{.,.,\}$ and $\{.,.\}^a$ linked by :

$$\{f,g\}^a= {1 \over a} \{af,ag\}$$

\noindent for any two $L$-admissible functions $f$ and $g$.
}
\end{thm}

\noindent{\sl Proof:}  The bracket $\{.,.\}$ is defined on 
 $A^{\infty}(L)$ by: 

$$\{f,g\}= - \lag e_f,e_g \rg_-,$$

\noindent where $e_f=(X_f, \varphi_f)+(df,f)$ and 
$e_g=(X_g, \varphi_g)+(dg,g)$ are sections of $L$. 
 By Remark~\ref{product}, the functions $af$ and $ag$ are 
 $L$-admissible. 
 Thus, the expression $\{af,ag\}$ makes sense. Moreover, Proposition~\ref{algebra} ensures
  that $\{.,.\}$ is a Lie bracket.
 Hence it remains to prove that $\{.,.\}^a$ satisfies the Jacobi identity. 
 For any $L$-admissible functions $f_1$, $f_2$ and $f_3$, we have

$$\{\{f_1,f_2\}^a, f_3\}^a ={1 \over a} \{a\{f_1,f_2\}^a,af_3\}
={1 \over a} \{\{af_1,af_2\},af_3 \}.$$

\noindent The Jacobi identity for the bracket $\{.,.\}$ implies

$$\{\{f_1,f_2\}^a, f_3\}^a+c.p.={1 \over a}\{\{af_1,af_2\},af_3 \}+c.p.=0.$$

\noindent This concludes the proof of Theorem~\ref{step2}.

\smallskip

\hfill \qed

 \medskip

  In fact, $\{.,.\}^a$ is related to a specific isotropic 
 sub-bundle of  $({\cal E}^1(M), \lag .,, \rg_+)$
  which depends on $a$. To define properly this sub-bundle,
 we need the following lemma:
 
\begin{lemma}
\label{step1}
 {\sl Let $L$ be a maximal isotropic sub-bundle
 of  $({\cal E}^1(M), \lag .,, \rg_+)$. Let $a$ be a smooth function on $M$ 
 that vanishes nowhere. Fix a vector field $Y$ and consider
 the sub-bundle $\widehat L \subset {\cal E}^1(M)$ 
 whose fibre at a point $x \in M$ is given by
$$\widehat L_x=span \Big\{ \big((a X +f Y)_x, (a \varphi-i_{Y} \alpha)_x \big) + \big( \alpha_x,f_x \big) \ | \ (X, \varphi)+(\alpha,f) \in \Gamma(L)
  \Big \}.$$

\noindent Then $\widehat L$ is  maximally isotropic. 
 }
\end{lemma}

\noindent{\sl Proof:} 
Consider two elements of $\Gamma (\widehat L)$ denoted by  
 $\delta_1$ and $\delta_2$ such that
$$\delta_n=(a X_n +f_n Y,  \ a \varphi_n-i_{Y}\alpha_n) + (\alpha_n, f_n),
\ \ \forall n=1,2$$
 
\noindent with $e_n=(X_n, \varphi_n) + (\alpha_n, f_n) \in \Gamma(L).$
 We have

$$ \lag \delta_1, \delta_2 \rg_+ = a \lag e_1, e_2 \rg_+.$$

\noindent Since $L$ is isotropic, we deduce that $ \lag \delta_1, \delta_2 \rg_+=0$. This proves that $\widehat L$ is isotropic under $ \lag .,. \rg_+$. This is a maximal isotropic sub-bundle of ${\cal E}^1(M)$. 
 Indeed, suppose that $L'$ is isotropic sub-bundle of ${\cal E}^1(M)$ which 
 contains  $\widehat L$. Then for any $\delta_1 \in \Gamma(\widehat L)$ and 
 for any 
 $e'=(Z,h)+(\beta, k) \in \Gamma(L')$, we have

\begin{eqnarray*}
 0=2 \lag \delta_1, \ e' \rg_+&=& i_{(a X_1 +f_1 Y)} \beta
  +i_{Z} \alpha_1 + hf_1 +k(a \varphi_1-i_{Y}\alpha_1) \cr
&= &  \lag e_1, \ (Z-kY, h+i_{Y} \beta)
 +(a\beta, ak)\rg_+,  
\end{eqnarray*}

\noindent where $e_1=(X_1, \varphi_1)+(\alpha_1, f_1)$ is the section of $L$ corresponding to $\delta_1$.
 Since $L$ is maximally isotropic, we deduce that 
$(Z-kY, h+i_{Y} \beta)+a(\beta, k)$ is in $\Gamma(L)$. Hence we obtain
the existence of a section $(X', \varphi')+(\alpha', f') \in \Gamma(L)$ such
 that
$$(Z-kY, h+i_{Y} \beta)+a(\beta, k)=a (X', \varphi')+a(\alpha', f').$$

\noindent This means that
$$e'=(aX'+f'Y,\ a \varphi' -i_{Y} \alpha'_i)+(\alpha',f').$$

\noindent Therefore  $\widehat L=L'$. This proves the lemma.

\hfill \qed

\begin{prop}
\label{bracket}
 {\sl Let $L$ be a maximally isotropic sub-bundle of  
 $({\cal E}^1(M), \lag .,, \rg_+)$ and let $a$ be a smooth function on $M$ 
 that vanishes nowhere. Assume that $(Y_a, \theta_a)+(da,0)$ is
 a smooth section of $L$. Consider the sub-bundle $L^a \subset {\cal E}^1(M)$ 
 whose fibre at a point $x \in M$ is given by
$$L^a_x=span \Big\{ \big((a X +f Y_a)_x, (a \varphi-i_{Y_a} \alpha)_x \big) + \big( \alpha_x,f_x \big) \ | \ (X, \varphi)+(\alpha,f) \in \Gamma(L) \Big\}.$$
\noindent Then the corresponding bracket on  $A^{\infty}(L^a)$ coincide
with $\{.,.\}^a$. 
 }
\end{prop}

\noindent{\sl Proof:} We have to prove that if 
 $\delta_f=(Y_f, \psi_f) + (df, f)$ and 
  $\delta_g=(Y_g , \psi_g) + (dg, g)$ are in $\Gamma(L^a)$, 
  then 

  $$\{f,g\}^a= - \lag\delta_f,\delta_g \rg_-.$$

 \noindent Obviously, a function $f$ is $L$-admissible if and only
 if it is $L^a$-admissible. More precisely,  any section
 $(X_f, \varphi_f)+(df, f) \in \Gamma(L)$  corresponds to
 a section $\delta_f$ of $L^a$ given by

 $$\delta_f=(X_{af}, \varphi_{af}-f\theta_a -Y_a \cdot f)+ (df, f).$$

\noindent where

 $$X_{af}=aX_f +fY_a, \quad  \quad  \varphi_{af}=
 a \varphi_f+f\theta_a $$

\noindent and $ e_{af}=(X_{af}, \ \varphi_{af})+(daf,af) \in \Gamma(L).$ 
 We have a similar expression for $\delta_g$. Therefore,

\begin{eqnarray*}
2 \lag \delta_f,\delta_g \rg_- &=X_{ag} \cdot f -X_{af} \cdot g +f\varphi_{ag}
-g\varphi_{af}+i_{Y_a} (gdf-fdg). \cr
 \end{eqnarray*}

\noindent This is equivalent to the following equation:
$$ 2 \lag \delta_f,\delta_g \rg_-={1 \over a}
 \Big(\lag e_{af}, e_{ag} \rg_+ \
 +i_{(gX_{af}-fX_{ag})}da\Big) +i_{Y_a}(gdf-fdg).$$
\noindent Now, we replace $X_{af}$ and $X_{ag}$ by their value
 and use the fact that 
   $$i_{Y_a}da= \lag (Y_a, \theta_a)+(da,0),  \ (Y_a, \theta_a)+(da,0) \rg_+=0.$$   
 \noindent Then we obtain 
$$\lag \delta_f,\delta_g \rg_-=-{1 \over a}\{af,ag\}
+ {1 \over 2} \Big(i_{(gX_{f}-fX_{g})}da +i_{Y_a}(gdf-fdg) \Big).$$

\noindent Since $fe_g-ge_f$ is in $\Gamma(L)$, we have

$$\lag fe_g-ge_f, \ (Y_a, \theta_a)+(da,0) \rg_+
=i_{(gX_{f}-fX_{g})}da +i_{Y_a}(gdf-fdg)=0.$$

\noindent We deduce that 
$$-\lag \delta_f,\delta_g \rg_-={1 \over a} \{af,ag \}=\{f,g\}^a.$$

\hfill \qed

 \medskip

 \noindent Now we are ready to state our second main result:

\begin{thm}
\label{Laconf}
 {\sl Let $L$ be a Dirac structure of ${\cal E}^1(M)$.
Assume that smooth sections of $L$ are locally spanned
  by  elements of the form $(X_f, \varphi_f)+(df, f)$.
 Let $a$ be a smooth function on $M$ that vanishes nowhere and such that
 an element of the type $(Y_a, \theta_a)+(da,0)$ is in $\Gamma(L)$.
  Let $L^a$ denote the sub-bundle whose fibre at a point $x \in M$ is given by:
$$L^a_x=span \Big\{ \big((a X +f Y_a)_x, (a \varphi-i_{Y_a} \alpha)_x \big) + \big( \alpha_x,f_x \big) \ | \ (X, \varphi)+(\alpha,f) \in \Gamma(L) \Big\}.$$

 \noindent Then, $L^a$ is a Dirac structure of ${\cal E}^1(M)$.
 }
\end{thm}

 \noindent To prove Theorem~\ref{Laconf}, we use the following lemma:

\begin{lemma}
\label{step3}
 {\sl Let $L$ be a  maximally isotropic sub-bundle of  ${\cal E}^1(M)$ with
  respect to $ \lag .,. \rg_+$. Then $L$ is a Dirac structure if and only if 
for any $e_1$, $e_2$ and $e_3 \in  \Gamma(L)$, we have:

$$ \lag \ [e_1,e_2], e_3 \rg_+=0.$$
 }
\end{lemma}

\noindent {\sl Proof:} This lemma is an immediate consequence of the fact that 
 $L$ is a maximally isotropic sub-bundle of  ${\cal E}^1(M)$.

\smallskip

\hfill \qed

\noindent {\sl Proof of Theorem~\ref{Laconf}:} 
 It follows from Lemma~\ref{step1} that $L^a$ is a maximal isotropic 
sub-bundle under $\lag .,. \rg_+$. To prove that $\Gamma(L^a)$ is 
closed under the extended Courant bracket $[.,.]$, we use Lemma~\ref{step3}.
 In fact, we can restrict ourselves to an open set ${\cal U}$ of $M$ 
 such that sections of $L$ over ${\cal U}$ are spanned by the family
 of sections $e_f=(X_f,\varphi_{f})+(df,f)$, where $f\in C^{\infty}({\cal U})$.  Therefore,  the set of smooth sections of $L^a$ over ${\cal U}$ 
  is spanned by elements of the form $ \delta_f=(Y_f, \psi_f)+(df,f)$, where 
$$Y_f=a X_f+f Y_a, \ \ \ \psi_f= a \varphi_{f}- Y_a \cdot f$$
 
 \noindent and $e_f=(X_f,\varphi_{f})+(df,f)$ is a smooth section of $L$.
 So, it is sufficient to show that 

$$ \lag \ [ \delta_{f_1},  \delta_{f_2}],  \delta_{f_3} \rg_+=0$$

\noindent for $L^a$-admissible functions $f_1$, $f_2$ and $f_3$.
Using Proposition~\ref{bracket}, we obtain  

$$ \lag \ [ \delta_{f_1},  \delta_{f_2}],  \delta_{f_3} \rg_+=
 \{\{f_1,f_2\}^a,f_3\}^a +c.p., $$

\noindent where $\{f_i,f_j\}^a= - \lag \delta_{f_i}, \delta_{f_j} \rg$.
 According to Theorem~\ref{step2}, $\{.,.\}^a$ is a Lie bracket on 
 $A^{\infty}(L)$ . Thus, the result is proved 

\smallskip

\hfill \qed

\begin{rmk} {\rm  In Theorem~\ref{Laconf}, if
  $\Gamma(L)$ is not locally spanned by the elements of the form
  $ e_f=(X_f, \varphi_f)+(df,f)$, one gets a maximal
 isotropic sub-bundle $L^a$ (see Lemma~\ref{step1}). But, it is not
 clear whether the set of sections of $L^a$ is closed under the bracket $[.,.]$ of ${\cal E}^1(M)$ or not. Observe that
    when $L$ corresponds to a Jacobi structure, 
  $\Gamma(L)$ is spanned by the elements 
 of the form $e_f$. We shall give another example where 
 this hypothesis of Theorem~\ref{Laconf} is fulfilled.
}
\end{rmk}
\medskip

\noindent {\bf Example.}  
 Assume that $(\omega,\eta)$ defines a locally conformal pre-symplectic 
 structure, where the rank of $\omega$ is constant. Denote 
 $D=\{  X \in TM \ | \ i_X \omega=0\}$ and $D^{\perp}$ its annihilator
 in $T^*M$. In fact, we have $D^{\perp}=\{i_X \omega \ | \ X \in TM\}$.
   Since $d\omega= \eta \wedge \omega$,  the distribution $D$ is integrable.
 %(by the classical theorem of Frobenius). 
 Hence, for any point $x$ of $M$, there exist coordinates 
 $(y_1, \dots , y_n)$ defined in a open neighborhood 
 ${\cal V}$ of $x$ such that $D^{\perp}|_{\cal V}$ is spanned by the 1-forms
  $dy_1, \dots dy_r$. \par 
 We denote by $L_{\omega,\eta}$ the Dirac structure associated with 
 $(\omega,\eta)$. So, we have
  $$L_{\omega,\eta}(x)=\{ (X, -i_X \eta)_x+(i_X \omega + f \eta, \ f)_x  \ |  \ X \in \Gamma(TM), \ f \in C^{\infty}(M) \}.$$
 
\noindent  From above, we deduce the existence of $r$ vector fields
 $X_1, \dots , X_r$ defined on ${\cal V}$ such that
 the set of smooth sections of $L_{\omega,\eta}$ over ${\cal V}$ is spanned by 
  $$ (0,0)+(\eta,1) \ \ {\rm and} \ \
  (X_k,-i_{X_k}\eta)+(dy_k,0), \ \ k=1,...r. $$  

\noindent We can suppose without loss of generality that 
 $\eta={dh \over h}$, where $h \in C^{\infty}({\cal V})$  
 and $h(x) > 0,$ for any $x \in {\cal V}$. Consider
  
 $$e_k=(X_k,-i_{X_k} \eta)+(dy_k+y_k \eta, y_k).$$

\noindent This  can be written as follows: 
$$e_k={1 \over h} \Big((Z_k,-i_{Z_k}\eta)+(dz_k, z_k),\Big), \ \forall k=1,...r,$$ 
 
 \noindent where

$$Z_k=hX_k \ \ \  {\rm and} \ \ \ z_k =h y_k, \ \forall k=1,...r.$$

\noindent  So, the set of smooth sections of $L_{\omega,\eta}$ over
 ${\cal V}$ is spanned by $(0,0)+(dh,h)$ and the $(Z_k,-i_{Z_k})+(dz_k, z_k),$ $\forall k=1,...r.$ We deduce that smooth sections of $L_{\omega,\eta}$
  are locally spanned  by  elements of the form $(X_f, \varphi_f)+(df, f)$. 

\hfill \qed

 Theorem~\ref{Laconf} gives rise to the following two definitions:

\begin{dfn}
\label{conf}
{ \rm
Let $a$ be a smooth function on $M$ which vanishes nowhere.
 Let $L$ and $L'$ be two Dirac structures of ${\cal E}^1(M)$. Assume that 
 there exists a smooth section  of the form $(Y_a, \theta_a)+(da,0)$ in
 $\Gamma(L)$  and the smooth sections of $L'$ are of the type
 $(a X+f Y_a, a \varphi-i_{Y_a} \alpha)+(\alpha, f)$, where
 $(X,\varphi)+(\alpha,f)$ belongs to $\Gamma(L)$.
 Then, $L'$ is said to be {\sl $a$-conformal} to $L$.

}
\end{dfn}

This first definition suggests to define a relation ${\cal R}$ among Dirac structures of
  ${\cal E}^1(M)$ by: $L' {\cal R} L$ if and only if there exists a smooth 
 function  $a$  satisfying  $a(x) \ne 0$, $\forall x \in M$ and $L'$ is 
 $a$-conformal to $L$. In such a case, we say that 
 $L'$ is {\sl conformally equivalent} to $L$ on $M$.\par

 \begin{dfn}{\rm
 A {\sl  conformal Dirac structure on $M$} is defined by an open cover 
 $\{U_i, i \in I\}$ of $M$  and a collection of Dirac structures 
 $L_i \subset {\cal E}^1(U_i)$ such that if $U_i \cap U_j$ is not empty
 then  $L_i$ is conformally equivalent $L_j$ upon $U_i \cap U_j$, where 
 ${\cal E}^1(U_i)= (TU_i \times \reals) \oplus  (T^*U_i \times \reals)$. \par

 Consider another open cover $\{V_r, r \in R\}$ of $M$
  such that each $V_r$ is endowed with a Dirac structure 
 $\widehat L_r \subset {\cal E}^1(V_r)$  and 
 whenever $V_r \cap V_s \ne  \emptyset$,
 $\widehat L_r$ is conformally equivalent to $\widehat L_s$.
  We say that this datum and $\{(U_i,L_i),  i \in I\}$  define the 
 same conformal Dirac structure on $M$  if  we have
 $$U_i \cap V_r \ne  \emptyset  \ \Longrightarrow \  L_i\ \mbox{\rm is 
 conformally equivalent to} \ \widehat L_r, \ \forall \ (i,r)\in I \times R .$$

}
\end{dfn}

 \noindent We have:

\begin{prop}{\sl   The relation ${\cal R}$ is an equivalence relation.
}
\end{prop}

 \noindent {\sl Proof:} It is clear that $L$ is $1$-conformal to itself,
where 1 is the constant function equal to 1 at any point of $M$.
 Now, assume that $L'$ is $a$-conformal to $L$. By definition,
there exists a smooth section $(Y_a, \theta_a)+(da,0) \in
 \Gamma(L)$. Since $Y_a \cdot a=0$ and any smooth
section of $L'$ can be obtained from an element of $\Gamma(L)$, we deduce
that $(aY_a, a \theta_a)+(da,0)$ belongs to $\Gamma(L')$. This implies
that $(Z_a, \eta_a)+(d({1 \over a}),0) \in  \Gamma(L')$, where 

$$Z_a=-{1 \over a}Y_a \quad {\rm and}
 \quad \eta_a=-{1 \over a}\theta_a.$$

\noindent Moreover we have

$$(X',\varphi')+(\alpha,f) \in \Gamma(L') \iff 
({1 \over a} X+f Z_a, \ {1 \over a} \varphi-i_{Z_a}  \alpha)+(\alpha, f) \in \Gamma(L).$$
\noindent This proves that $L$ is ${1 \over a}$-conformal to $L'$.
Furthermore, it is easy to see that if $L''$ is $b$-conformal to $L'$ and $L'$ is $a$-conformal to $L$, then $L''$ is $ab$-conformal to $L$. So, we obtain the
proposition.

\hfill \qed

\begin{rmk}
{\rm
 Obviously the equivalence class (with respect to ${\cal R}$)
  of any Dirac structure 
 $L \subset {\cal E}^1(M)$ defines a conformal Dirac structure on $M$. But
   a conformal Dirac structure may not be globally defined on a manifold $M$.   An interesting case is that of a conformal Dirac structure defined
  locally  by conformal pre-symplectic structures but not defined 
 globally. This case will be analyzed in a coming paper.  \par 

}
\end{rmk}

\noindent {\bf Acknowledgements.} I wish to thank A. Weinstein for 
  calling my attention to this problem and comments on an earlier version
 of this paper. Part of this work was conceived during a visit in the 
  Mathematics Department at UNC-Chapel Hill. I would like to express 
 my gratitude to this institution. I thank I. Assani for his hospitality. 
 Thanks are due  to the Abdus Salam International Centre for 
 Theoretical Physics for its support, also to J.-P. Dufour and J. Stasheff for
  discussions.

\bigskip

\noindent {\bf Current address:} \par
\noindent {\small DEPARTMENT OF MATHEMATICS,} 
 {\small UNIVERSITY OF NORTH CAROLINA,\par
  \noindent CHAPEL HILL, NC} 27599-3250.\par
  E-mail: aissaw@math.unc.edu

\end{document}